\documentclass[runningheads]{llncs}
\usepackage{amsmath,amssymb,amsthm,mathtools}
\usepackage{xcolor} %per usare \textcolor
\usepackage{comment}
\usepackage{hyperref}
\usepackage{graphicx, graphics, subcaption}
\usepackage[leftcaption]{sidecap}
\usepackage{booktabs,siunitx} % better tables
\usepackage{afterpage} % to put figure in the next page

\usepackage[font=footnotesize]{caption}
\usepackage[misc]{ifsym}
\newcommand{\corrAuthor}{$^{\left(\textrm{\Letter}\right)}$}

\usepackage{algorithm}% http://ctan.org/pkg/algorithms
\usepackage{algorithmic}% http://ctan.org/pkg/algorithms
\makeatletter
\renewcommand{\ALG@name}{Algorithm}
\renewcommand{\listalgorithmname}{List of \ALG@name s}
\makeatother
\usepackage{tikz}
\tikzset{mathterm/.style={draw=black,fill=white,rectangle,anchor=base}}
\tikzstyle{every picture}+=[remember picture]
\usepackage{pgfplots}
\usepackage{verbatim}
\theoremstyle{plain}
\newtheorem*{assumption*}{Assumption}

\newcommand{\Ucal}{\mathcal{U}}
\newcommand{\rn}{\mathbb{R}^n}
\newcommand{\rnn}{\mathbb{R}^{n \times n}}
\newcommand{\R}[2]{\mathbb{R}^{#1 \times #2}}

\newcommand{\MATLAB}{\textsc{Matlab} }
\renewcommand{\phi}{\varphi}

\begin{document}
\title{A New Algorithm for the LQR  Problem with Partially Unknown Dynamics\thanks{This research has been partially supported by the INdAM Research group GNCS.}}
\titlerunning{A New Algorithm for LQR Problem with Partially Unknown Dynamics}
\author{Agnese Pacifico\corrAuthor \and
Andrea Pesare \and
Maurizio Falcone}
\authorrunning{A. Pacifico et al.}
\institute{Department of Mathematics, Sapienza University of Rome, Rome, Italy\\
\email{pacifico.1699761@studenti.uniroma1.it, \{pesare,falcone\}@mat.uniroma1.it}}
\maketitle
\begin{abstract}
We consider an LQR optimal control problem with partially unknown dynamics. We propose a new model-based online algorithm to obtain an approximation of the dynamics {\em and} the control at the same time during a single simulation.
\end{abstract}
\section{Introduction}
The Linear-Quadratic Regulator (LQR) optimal control problem \cite{anderson2007LQR,Fleming-Rishel} is a classical problem in control theory with a wide range of applications. When the dynamics is known, the optimal control is obtained in feedback form solving a backward Riccati differential equation.

Some Reinforcement Learning (RL) problems can be seen as LQR problems where the dynamics is partially or completely unknown. Some RL algorithms try to learn a model for the dynamics and use this model to find the optimal policy. These are called \textit{model-based algorithms}. Others recover the optimal control directly, without reconstructing the dynamics, and these are called \textit{model-free algorithms}. For a broad overview on RL, we recommend \cite{SuttonBarto}.

The connection between RL and optimal control theory was already identified in the past \cite{sutton1992}. Recently, Palladino and co-authors tried to clarify this relationship via some rigorous proofs, identifying some RL tasks as real optimal control problems with unknown dynamics \cite{murray2018model,pesare2020LQR,pesare2021ECC}. In this context, we propose a {\em new model-based algorithm for LQR problems} where the dynamics is partially unknown, which takes contributions from both fields. In fact, our algorithm can be considered as a case of Bayesian RL, a class of model-based algorithms where the controller builds a stochastic model of the dynamics and updates it according to Bayesian statistics \cite{BRL-Survey}. On the other hand, we borrow the LQR solution from optimal control theory, to get the synthesis of a suitable control \cite{anderson2007LQR}. 

In particular, our algorithm is similar to PILCO \cite{PILCO}, from which it takes the use of Gaussian distributions and the whole process of Bayesian update. However, we propose here some novelties. The first is that in PILCO the optimal control is chosen through a minimization procedure in a class of controls, whereas we solve a Riccati differential equation. The second is that  PILCO needs several trials to reconstruct the dynamics and stabilise the system. Our algorithm is designed to approximate the dynamics \textit{and} to find a suitable control in a single trial but can be applied only to linear dynamical systems. Finally, let us mention that other works already dealt with LQR problems with unknown dynamics (see i.e. \cite{fong2018dual,frueh2000LQL,li2019lq,pang2019adaptive} and references therein), but they all need several trials to converge, whereas our method works with just one simulation.

The paper is structured as follows. In Section 2 we recall the LQR problem. In Section 3 we present our algorithm and discuss some implementation details. Finally, in Section 4 we show and discuss some numerical tests.

\section{The classical LQR problem}\label{sec:classical_LQ}
The Linear-Quadratic Regulator (LQR) problem \cite{anderson2007LQR,Fleming-Rishel} is an optimal control problem with linear dynamics and quadratic cost. In the finite horizon case, the state of the system $x(t) \in \rn$ evolves according to the following controlled dynamics
\begin{equation}\label{eq:problem_A}
\left\{ 
\begin{array}{lll}
\dot{x}(t)&=&\widehat{A} x(t) + Bu(t),\quad t \in [0,T] \\
x(0)&=&x_0 .
\end{array}
\right.
\end{equation}
We will denote by $\R{m}{n}$ the space of matrices with $m$ rows and $n$ columns; $\mathbb{I}_n$ and $\mathbf{0}_n$ will be respectively the identity and zero matrices in $\rnn$. Here 
$\widehat{A} \in \R{n}{n}$ and $B\in \R{n}{m}$.
The control function $u(t) \in \mathbb{R}^m$ must be chosen among the admissible controls $\Ucal:=\{u:[0,T] \rightarrow \mathbb{R}^m\; \mathrm{Lebesgue}\; \mathrm{measurable}\}$ to minimize the quadratic cost functional
\begin{equation}\label{eq:cost_functional_A}
 J_{x_0}[u] \coloneqq \frac{1}{2} \left(\int_0^T \left( x(t)^T Q x(t) + u(t)^T R u(t) \right) dt + x(T)^T Q_f x(T) \right) ,
\end{equation}
where $Q, Q_f \in \rnn$ are symmetric and positive semi-definite, and $R \in \R{m}{m}$ is symmetric and positive definite.
For any $x_0 \in \mathbb{R}^n$, 
%the value function in $x_0$ is defined as \agn{$ V(x_0) = \inf_{u \in \Ucal} J_{x_0}[u] $,}
%\begin{equation}\label{eq:value_function_A}
%    V(x_0) = \inf_{u \in \Ucal} J_{x_0}[u],
%\end{equation}
we will call $\bar{u}$ an optimal control starting from $x_0$ if and only if 
\begin{equation}\label{eq:optimal_control_A}
    J_{x_0}[\bar{u}] \leq J_{x_0}[u] \quad \forall \, u \in \Ucal .
\end{equation}

When the dynamics is fully known, the optimal control can be obtained in feedback form \cite{anderson2007LQR}. Indeed, if  $P(t)$ is the unique symmetric solution of the Riccati differential equation
\begin{equation}
\left\{ 
\begin{aligned}
    -\dot{P}(t)&=\widehat{A}^TP(t)+P(t)\widehat{A}-P(t)BR^{-1}B^TP(t)+ Q, \quad t \in [0, T]\\
    P(T)&=Q_f ,
\end{aligned}
\right.
\end{equation}
the optimal control is given by $\bar{u}(t) = -R^{-1} B^T P(t) x(t)$.

We want to investigate: what happens if the dynamics is partially unknown? How can one find a suitable control? We considered the following framework:

\begin{assumption*} We assume that the matrices $B \in \R{n}{m}$, $Q, Q_f \in \R{n}{n}$ and \mbox{$R \in \R{m}{m}$} are given, whereas the state matrix $\widehat{A} \in \R{n}{n}$ is unknown. 
\end{assumption*}

\section{An Online Algorithm for the LQR Problem}\label{sec:algorithm}
In this section, we describe our model-based algorithm able to solve the LQR problem without knowing the matrix $\widehat{A}$. The algorithm is online, meaning that it doesn't need any previous simulation or computation. Our goal is twofold: first, we look for a good estimate for the unknown dynamics matrix $\widehat{A}$; and secondly, we want to choose a control that can stabilise the system. Furthermore, we assume that we can run the experiment just once, so the system must be controlled while the dynamics is still uncertain. To this end, we use a technique to get an estimate of the matrix $\widehat{A}$ and to update the control at the same time.

We divide the interval $[0,T]$ into equal time steps of length $\Delta t$, globally we will have $N$ time steps  and we group them into rounds of $S$ steps each. The $i$-th round will be denoted by $[t_{i-1},t_i]$ and a superscript index will indicate a single time step:
\[ t_i^j = t_{i-1} + j \Delta t , \quad j=0, \dots, S . \]

The current knowledge of the dynamics matrix is represented as a probability distribution over matrices. For each round, two major operations are carried out:
\begin{enumerate}
    \item At the beginning of each round, a probability $\pi_{i-1}$ is given from the previous round. We use the mean $\bar{A}_{i-1}$ of this distribution and compute a feedback control for the round solving a Riccati equation;
    \item At the end of each round, the current probability distribution is updated using Bayesian formulas, according to the trajectory observed during the round. The output is a new probability distribution $\pi_i$.
\end{enumerate}

The whole algorithm is summarised below as Algorithm~\ref{alg:RLonline}. In the following, we will give more technical details.
\afterpage{
\begin{algorithm}[t!]
\caption{An online algorithm for the LQR problem}
\label{alg:RLonline}  
\begin{algorithmic} 
    \STATE{Divide $[0,T]$ into $M$ intervals of length $\Delta t$, with $\Delta t = \frac{T}{M}$;}
    \STATE{Group the intervals in \emph{rounds}, each containing $S$ intervals;}
    \STATE{Choose a prior distribution $\pi_0$ over matrices;}
    \FOR{$i$ from 1 to the \emph{rounds} number}
        \STATE{Find a feedback control $u^*_i$ solving a Riccati eq. with the mean matrix $\bar{A}_{i-1}$;}
        \STATE{Use $u^*_i$ as control for all the steps in the $i$-th \emph{round};}
        \STATE{Observe the actual trajectory;}
        \STATE{Update $\pi_i$ according to the data from the observed trajectory}
    \ENDFOR  
\end{algorithmic}  
\end{algorithm}}

\subsubsection*{Prior distribution}
The algorithm requires the choice of a prior distribution $\pi_0$. We fix some $m_0\in\rn$ and $\Sigma_0\in\rnn$ and consider a random matrix $A$ such that each of its rows $r_k$ is distributed as an independent gaussian vector with mean $m_0$ and covariance matrix $\Sigma_0$
$$A=
\begin{bmatrix}
 & r_1^T & \\
 & \vdots & \\
 & r_n^T & \\
\end{bmatrix}
\qquad \qquad r_k\sim\mathcal{N}(m_0,\Sigma_0) \qquad \forall \, k = 1, \dots, n$$ 
A typical choice for the prior distribution, when no information about the true matrix $\widehat{A}$ is available, is $m_0 = \left( 0, \dots, 0 \right)^T$ and $\Sigma_0 = n \, m \, \mathbb{I}_n$.

\subsubsection*{Feedback control}
At the beginning of the round $[t_{i-1},t_i]$ our knowledge of the matrix $\widehat{A}$ is described by the distribution $\pi_{i-1}$. In order to find the control to apply, we solve the evolutive \emph{Riccati} equation associated with the matrix $\bar{A}_{i-1}$, where $\bar{A}_{i-1}$ is the mean of the distribution $\pi_{i-1}$. The Riccati equation reads
\begin{equation}\label{eq:LQRiccati}
\left\{ 
\begin{aligned}
    -\dot{P}(t)&=\bar{A}_{i-1}^TP(t)+P(t)\bar{A}_{i-1}-P(t)BR^{-1}B^TP(t)+ Q \quad t \in [t_{i-1}, T]\\
    P(T)&=Q_f .
\end{aligned}
\right.
\end{equation}
If we denote by $P_i(t)$ the solution of \eqref{eq:LQRiccati}, our control will be the feedback control given by
\begin{equation}
    u_i^*(t_i^j)=-R^{-1}B^T P_i(t_i^j)x^j_i \qquad j=0,\ldots S-1 ,
\end{equation}
where $x_i^j=x(t_i^j)$ is the observed state of the system. Since the control must be defined for all instants $t \in [t_{i-1}, t_i]$, we choose a piecewise constant control: $u_i^*(t) = u_i^*(t_i^j)$ for $t \in [t_i^j, t_i^{j+1}]$.
Note that we need real-time observations of the system state to compute $u^*_i$, since it depends on $x$. 

\subsubsection*{Distribution update}
At the end of each round we can update the probability distribution using {\em Bayesian regression formulas}. More precisely, we know that during the time step $[t_i^j, t_i^{j+1}]$ the system evolves according to \eqref{eq:problem_A} where we plug-in the chosen piecewise constant control.
We can easily get an approximation of the state derivative with a finite difference scheme using the state observations $x_i^j=x(t_i^j)$. By a first-order scheme, rearranging the terms,  we get
\begin{equation}\label{eq:regression_y}
    \widehat{A}x_i^j \simeq \frac{x_i^{j+1}-x_i^{j}}{\Delta t}-Bu^*(t_i^j) 
    \qquad j=0,\ldots,S-1 .
\end{equation}
We interpret these data as observations of the dynamics with a Gaussian noise due to the error in the derivative approximation and the measurements. We treat each row $\widehat{r}_k$ of $\widehat{A}$ separately.
Denoting by $y^j$ the right hand side in \eqref{eq:regression_y}, we can write its $k$-th component as
$$y^j_{(k)}=\widehat{r}_k x_i^j+\varepsilon \qquad \text{with }\varepsilon\sim\mathcal{N}(0,\sigma^2) ,$$
where $\sigma$ is a parameter chosen according to the order of the derivative approximation. 
%We have a prior distribution for it given by $\pi_{i-1}$, so we can compute the posterior distribution $\pi_{i}$. 

For each row $k$ we have a prior distribution $r_k\sim\mathcal{N}(m_0,\Sigma_0)$ given by $\pi_{i-1}$. We define $X$ as the matrix whose columns are $x_i^0,\ldots,x_i^{S-1}$ and $y$ as the column vector  $y^0_{(k)},\ldots,y^{S-1}_{(k)}$. Therefore we obtain a posterior distribution \mbox{$r_k|_{y,X}\sim\mathcal{N}(m,\Sigma)$} where $\Sigma^{-1}=\frac{1}{\sigma^2}XX^T+\Sigma_0^{-1}$ and $m=\Sigma ( Xy/\sigma^2 +\Sigma_0^{-1}m_0 )$.
For more details about Bayesian Linear Regression see for instance the extensive monography by Rasmussen and Williams \cite{Rasmussen2006}.

\subsubsection*{Higher-order schemes}
%Note that the derivative in \eqref{eq:regression_y} can be approximated with higher order formulas, for example
%\[ \dot{x}(t_i^j)\simeq\frac{x^{j+1}-x^{j-1}}{2 \, \Delta t} \qquad \text{or} \qquad \dot{x}(t_i^j)\simeq\frac{-x^{j+2}+8x^{j+1}-8x^{j-1}+x^{j-2}}{12 \, \Delta t} \]
The derivative in \eqref{eq:regression_y} can be approximated with higher order finite difference schemes \cite{strikwerda2004}. Note that trajectory regularity is required in the interval containing the nodes used in the approximation. While first-order approximation uses only two nodes, higher-order approximations use more nodes, thus the control cannot jump at each time step and we have to keep it constant for more steps.

\begin{remark}[Heuristic argument for convergence]\label{rmk:heuristics}
The algorithm cannot find the optimal control for the problem, since at the beginning the matrix $\widehat{A}$ is unknown and it needs at least some steps to have a good estimate for $\widehat{A}$. However, from Bayesian Regression theory (see \cite{Rasmussen2006}) we know that the more data we observe, the more precise our distribution $\pi_i$ becomes, eventually converging to the Dirac delta $\delta_{\widehat{A}}$.

Furthermore, let us note that since $\pi_i$ is converging to $\delta_{\widehat{A}}$, then also the mean matrix $\bar{A}_i$ of $\pi_i$ is converging to $\widehat{A}$. Thus, after few rounds, the feedback control computed by the algorithm (which is using $\bar{A}_i$) in the interval $[t_i, T]$ should be close to the optimal control of a trajectory of the real dynamics which starts from the same point $x_{i-1}^S$.

Finally, when $\Delta t \to 0$, the algorithm reaches earlier a good estimate of the dynamics matrix. This means that also the computed control is closer to the optimal one. All these heuristics are confirmed by the numerical simulations in the next section.

\end{remark}

\section{Some Numerical Tests}\label{sec:numerical_tests}
The following numerical tests for the algorithm described above were performed in \MATLAB and took few seconds to run.
\subsubsection*{Test 1}
We first consider a dynamical system where the state lies in $\mathbb{R}^2$ and the control is 1-dimensional, i.e. $n=2$ and $m=1$. The LQR problem is defined by the following matrices:
$$\widehat{A}=
\begin{pmatrix*}[r]
0 & \phantom{-}1 & \\
-1 & \phantom{-}0 &
\end{pmatrix*}
\quad
B=\begin{pmatrix*}
& 0 &\\
& 1 &
\end{pmatrix*}
\quad
Q=
\begin{pmatrix*}[r]
& 1 & \phantom{1}0 &\\
& 0 & 1 &
\end{pmatrix*}
\quad
R=0.1
\quad
Q_f=\begin{pmatrix}
& 0 & \phantom{1}0 &\\
& 0 & \phantom{1}0 &
\end{pmatrix} .
$$
The time horizon is set to $T=5$ and the starting point is $x_0=(0, 1)^T$.
We assume that the matrix $\widehat{A}$ is unknown to the algorithm, though we use it to simulate the dynamics. We choose the prior distribution as recommended in Section~\ref{sec:algorithm}, using $m_0 = (0,0)^T$ and $\Sigma_0 = 2 I_2$; for all the tests we set $\sigma=\sqrt{10\Delta t^p}$ and $S=2p$, where $p$ is the order of the scheme used in the derivative approximation.

Fig.~\ref{fig:2x2ordine1} shows the \textit{piecewise constant} controls chosen by the algorithm with $p=1$ for different values of $\Delta t$ and the corresponding trajectories. %Notice that the controls are taken piecewise constant. 
Recall that the matrix $\widehat{A}$ is completely unknown at the beginning, so the control we apply in the first steps depends only on the prior distribution we have chosen and clearly is not accurate. This causes the trajectory to deviate from the optimal solution. However, after few steps the matrix $\widehat{A}$ is well approximated and the algorithm manages to stabilise the system anyway, bringing the state towards the origin. In Table~\ref{tab:errors_1} we have reported the cost of the solution found by the algorithm for different choices of $\Delta t$. When $\Delta t$ is smaller, the algorithm recovers the matrix $\widehat{A}$ quickly and thus the deviation from the optimal solution is smaller. This confirms the heuristics of Remark~\ref{rmk:heuristics}. We can also observe a numerical order of convergence equal to~1.
\begin{figure}[tbhp]
\begin{subfigure}[t]{0.32\textwidth}
\raisebox{-\height}{
\begin{tikzpicture}[scale=0.5]
	\begin{axis}[legend style={anchor=south, at={(0.6,0.05)}}]
		\addplot[color=blue] 	table{immagini/simulazione2x21_ord1_dt01/simulazione1_u_tru1.txt};
		\addplot[color=red]	table{immagini/simulazione2x21_ord1_dt01/simulazione1_u_unc_ru1.txt};
		\legend{optimal control, our control}
	\end{axis}
\end{tikzpicture}
}
\vspace{-7pt}
\caption{control $\Delta t =0.1$}
%\label{fig:controllo_ord1}
\end{subfigure}
\begin{subfigure}[t]{0.32\textwidth}
\raisebox{-\height}{
\begin{tikzpicture}[scale=0.5]
	\begin{axis}[legend style={anchor=south, at={(0.6,0.05)}}]
		\addplot[color=blue] 	table{immagini/simulazione2x21_ord1_dt005/simulazione1_u_tru1.txt};
		\addplot[color=red]	table{immagini/simulazione2x21_ord1_dt005/simulazione1_u_unc_ru1.txt};
		\legend{optimal control, our control}
	\end{axis}
\end{tikzpicture}
}
\vspace{-7pt}
\caption{control $\Delta t =0.05$}
\end{subfigure}
\begin{subfigure}[t]{0.32\textwidth}
\raisebox{-\height}{
\begin{tikzpicture}[scale=0.5]
	\begin{axis}[legend style={anchor=south, at={(0.6,0.05)}}]
		\addplot[color=blue] 	table{immagini/simulazione2x21_ord1_dt0025/simulazione1_u_tru1.txt};
		\addplot[color=red]	table{immagini/simulazione2x21_ord1_dt0025/simulazione1_u_unc_ru1.txt};
		\legend{optimal control, our control}
	\end{axis}
\end{tikzpicture}
}
\vspace{-7pt}
\caption{control $\Delta t =0.025$}
\end{subfigure}

\vspace{9pt}

\begin{subfigure}[t]{0.32\textwidth}
\raisebox{-\height}{
\begin{tikzpicture}[scale=0.5]
	\begin{axis}
		\addplot[color=blue] 	table{immagini/simulazione2x21_ord1_dt01/simulazione1_x_tru12.txt};
		\addplot[color=red]	table{immagini/simulazione2x21_ord1_dt01/simulazione1_x_unc12.txt};
		\node at (axis cs:0,1) {\large\textcolor{red}{$\bullet$}};
		\legend{optimal trajectory, our trajectory}
	\end{axis}
\end{tikzpicture}
}
\vspace{-7pt}
\caption{trajectory $\Delta t =0.1$}
\end{subfigure}
\begin{subfigure}[t]{0.32\textwidth}
\raisebox{-\height}{
\begin{tikzpicture}[scale=0.5]
	\begin{axis}
		\addplot[color=blue] 	table{immagini/simulazione2x21_ord1_dt005/simulazione1_x_tru12.txt};
		\addplot[color=red]	table{immagini/simulazione2x21_ord1_dt005/simulazione1_x_unc12.txt};
		\node at (axis cs:0,1) {\large\textcolor{red}{$\bullet$}};
		\legend{optimal trajectory, our trajectory}
	\end{axis}
\end{tikzpicture}
}
\vspace{-7pt}
\caption{trajectory $\Delta t =0.05$}
\end{subfigure}
\begin{subfigure}[t]{0.32\textwidth}
\raisebox{-\height}{
\begin{tikzpicture}[scale=0.5]
	\begin{axis}
		\addplot[color=blue] 	table{immagini/simulazione2x21_ord1_dt0025/simulazione1_x_tru12.txt};
		\addplot[color=red]	table{immagini/simulazione2x21_ord1_dt0025/simulazione1_x_unc12.txt};
		\node at (axis cs:0,1) {\large\textcolor{red}{$\bullet$}};
		\legend{optimal trajectory, our trajectory}
	\end{axis}
\end{tikzpicture}
}
\vspace{-7pt}
\caption{trajectory $\Delta t =0.025$}
\end{subfigure}
\vspace{5pt}
\caption{Simulations of Test 1 for different values of $\Delta t$. The first row shows the control chosen by the algorithm as a function of time (in red) compared with the optimal control (in blue), computed knowing the matrix $\widehat{A}$; the second row shows the trajectories in $\mathbb{R}^2$. The red dot is the trajectory starting point.}
\label{fig:2x2ordine1}
\end{figure}
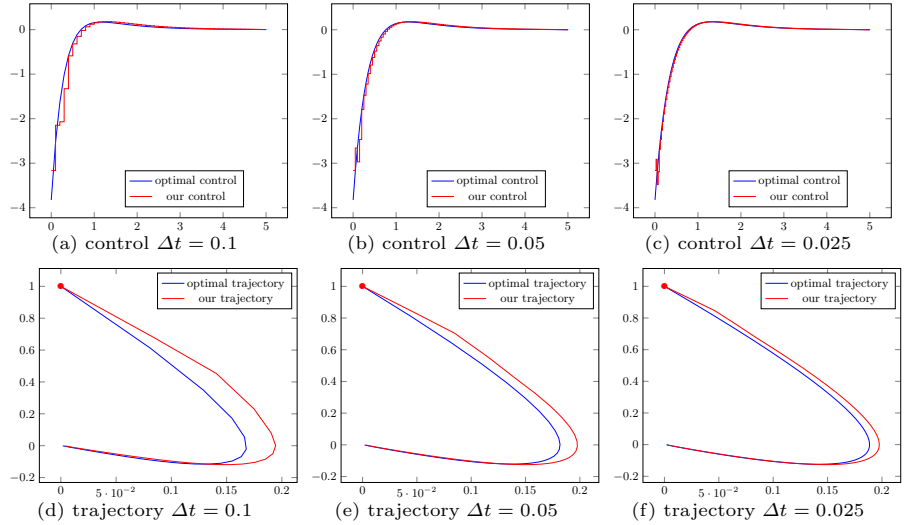
\begin{table}[tbhp]
    \caption{Numerical results for Test 1. (a): Cost of the trajectory for different values of $\Delta t$. The cost is compared with the optimal cost computed knowing $\widehat{A}$ ($C^*$ last row). (b): Error for $\widehat{A}$ after the simulation using different $\Delta t$ and different finite difference schemes for the derivatives approximation; $p$ is the scheme order.}
    \begin{subtable}[t]{.4\linewidth}
        \centering
        \caption{}
        \begin{tabular}{c S[table-format=2.4] S[table-format=2.4] c}
        \toprule
        \multicolumn{1}{c}{$\Delta t$} & \multicolumn{1}{c}{Cost} & \multicolumn{1}{c}{Error} & \multicolumn{1}{c}{Order}\\
        \midrule
        0.1            & 0.3897 & 0.0063 & - \\
        0.05           & 0.3866 & 0.0032 & 0.98 \\
        0.025          & 0.3849 & 0.0015 & 1.09 \\
        \midrule
        $C^*$           & 0.3834 & \multicolumn{1}{c}{-} & -\\
        \bottomrule
        \end{tabular}
        \label{tab:errors_1}
    \end{subtable}%
    \begin{subtable}[t]{.6\linewidth}
        \centering
        \caption{}
        \begin{tabular}{c|c c|c c|c c}
        \toprule
        & \multicolumn{2}{c|}{$p=1$} & \multicolumn{2}{c|}{$p=2$} & \multicolumn{2}{c}{$p=4$}\\
        \midrule
        $\Delta t$ & Error & Order & Error & Order & Error & Order \\
        \midrule
        0.1   & 0.167 & -	 & 9.74e-3 & -	 & 1.63e-4 & - \\
        0.05  & 0.087 & 0.94 & 1.91e-3 & 2.3 & 3.00e-6 & 5.8 \\
        0.025 & 0.045 & 0.95 & 4.14e-4 & 2.2 & 1.03e-7 & 4.8 \\
        0.01  & 0.018 & 1.00 & 6.30e-5 & 2.1 & 1.70e-9 & 4.5 \\
        \bottomrule
        \end{tabular}
        \label{tab:orders_24}
    \end{subtable}
\end{table}

We tried different finite difference schemes for the approximation of the state derivatives (see \textit{"Higher-order schemes"} in Section~\ref{sec:algorithm}). Table \ref{tab:orders_24} shows the error in the approximation of $\widehat{A}$ at the end of the simulation, when using schemes of order $p=1,2,4$ and for different values of $\Delta t$. As expected, when we consider more accurate approximations of the gradient, we get better estimations of the matrix $\widehat{A}$. Unfortunately, the same does not hold for the solution costs, which are not significantly improved if compared with the ones found by the first-order approximation.

\subsubsection*{Test 2}
\begin{figure}[tbhp]
\begin{subfigure}[t]{0.32\textwidth}
%\raisebox{-\height}{
\begin{tikzpicture}[scale=0.5]
	\begin{axis}
		\addplot[color=blue] 	table{immagini/simulazione4/simulazione1_u_tru1.txt};
		\addplot[color=red]	table{immagini/simulazione4/simulazione1_u_unc1.txt};
		\legend{optimal control, our control}
	\end{axis}
\end{tikzpicture}
%}
%\caption{prima componente}
\end{subfigure}
\begin{subfigure}[t]{0.32\textwidth}
%\raisebox{-\height}{
\begin{tikzpicture}[scale=0.5]
	\begin{axis}[legend style={anchor=south, at={(0.6,0.05)}}]
		\addplot[color=blue] 	table{immagini/simulazione4/simulazione1_u_tru2.txt};
		\addplot[color=red]	table{immagini/simulazione4/simulazione1_u_unc2.txt};
		\legend{ optimal control, our control}
	\end{axis}
\end{tikzpicture}
%}
%\caption{seconda componente}
\end{subfigure}
\begin{subfigure}[t]{0.32\textwidth}
%\raisebox{-\height}{
\begin{tikzpicture}[scale=0.5]
	\begin{axis}[legend style={anchor=south, at={(0.6,0.05)}}]
		\addplot[color=blue] 	table{immagini/simulazione4/simulazione1_u_tru3.txt};
		\addplot[color=red]	table{immagini/simulazione4/simulazione1_u_unc3.txt};
		\legend{ optimal control, our control}
	\end{axis}
\end{tikzpicture}
%}
%\caption{terza componente}
\end{subfigure}
\caption{The three components of the control in Test 2}
\vspace{-20pt}
\label{fig:43controlli}
\end{figure}
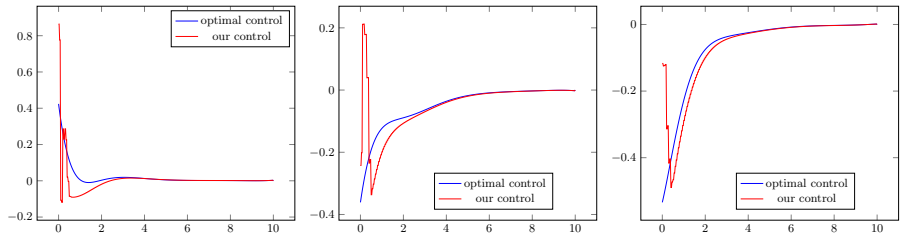

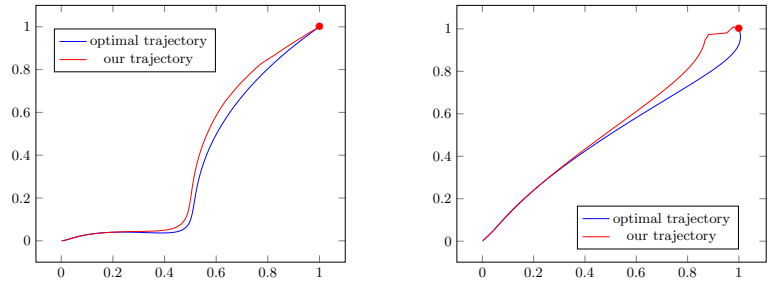
\begin{figure}[tbhp]
\hspace{25pt}
\begin{subfigure}[t]{0.49\textwidth}
%\raisebox{-\height}{
\begin{tikzpicture}[scale=0.6]
	\begin{axis}[legend style={anchor=south, at={(0.32,0.74)}}]
		\addplot[color=blue] 	table{immagini/simulazione4/simulazione1_x_tru34.txt};
		\addplot[color=red]	table{immagini/simulazione4/simulazione1_x_unc34.txt};
		\node at (axis cs:1,1) {\large\textcolor{red}{$\bullet$}};
		\legend{ optimal trajectory, our trajectory}
	\end{axis}
\end{tikzpicture}
%}
%\caption{Prime due componenti della traiettoria}
\end{subfigure}
\hspace{-15pt}
\begin{subfigure}[t]{0.49\textwidth}
%\raisebox{-\height}{
\begin{tikzpicture}[scale=0.6]
	\begin{axis}[legend style={anchor=south, at={(0.65,0.05)}}]
		\addplot[color=blue] 	table{immagini/simulazione4/simulazione1_x_tru12.txt};
		\addplot[color=red]	table{immagini/simulazione4/simulazione1_x_unc12.txt};
		\node at (axis cs:1,1) {\large\textcolor{red}{$\bullet$}};
		\legend{ optimal trajectory, our trajectory}
	\end{axis}
\end{tikzpicture}
%}
%\caption{Terza e quarta componente della traiettoria}
\end{subfigure}
\caption{The trajectory of Test 2 in $\mathbb{R}^4$, represented by projecting components in couples: $(x_1,x_2)$ on the left and $(x_3,x_4)$ on the right. The red dots indicate the starting point of the trajectory.}
\label{fig:43traiettorie}
\end{figure}

For the second test we choose $n=4$, $m=3$ and $T=10$. Our matrices are 
\[ \widehat{A}=\begin{pmatrix*}[r]
-0.0215 & -0.7776 & -0.1922 & 0.9123\\
-0.3246 & 0.5605 & -0.8071 & 0.1504\\
0.8001 & -0.2205 & -0.7360 & -0.8804\\
-0.2615 & -0.5166 & 0.8841 & -0.5304
\end{pmatrix*} \! , \
B=\begin{pmatrix*}[r]
-0.2937 & -0.6620 & -0.0982\\
0.6424 & 0.2982 & 0.0940\\
-0.9692 & 0.4634 & -0.4074\\
-0.9140 & 0.2955 & 0.4894
\end{pmatrix*} \! , \]
$Q=\frac{1}{4} \mathbb{I}_4$, $R=\frac{1}{3} \mathbb{I}_3$ and $Q_f=\mathbb{I}_4$, and the starting point is $x_0 = (1,1,1,1)^T$.
We set $\Delta t=0.025$, $S=4$ and use a second order approximation for the derivatives. 
Fig.~\ref{fig:43controlli} shows the control found by the algorithm and Fig.~\ref{fig:43traiettorie} the corresponding trajectory. The behaviour observed in Test 1 is even more visible here: in the first steps the control is not accurate since we do not know  $\widehat{A}$, but after few steps the algorithm learns more on the matrix and manages to stabilise the system. The cost of the control found by the algorithm is 1.111, whereas
the optimal cost computed using $\widehat{A}$ is 1.056.

\section{Conclusions}
We proposed a new algorithm designed to deal with LQR problems when the dynamics is partially unknown. Numerical tests presented in Section \ref{sec:numerical_tests} showed how it manages to approximate the dynamics \emph{and} to find a suitable control that brings the system towards the origin in a single simulation. Future works include the convergence analysis of the algorithm and possible extensions to the nonlinear case.

 \bibliographystyle{splncs04}
 \bibliography{LSSC21.bib}
\end{document}